\documentclass[a4paper,11pt]{article}
\usepackage[utf8]{inputenc}
\usepackage[T1]{fontenc}
\usepackage[dvipdfmx]{graphicx} 
\usepackage{amsmath,amsthm,amssymb}
\usepackage{geometry}
\usepackage{url}
\usepackage{hyperref}
\usepackage{float}
\usepackage{booktabs}
\usepackage{multirow}
\usepackage[numbers]{natbib}
\usepackage[x11names]{xcolor}
\usepackage{textcomp}
\usepackage{algorithm}
\usepackage{algpseudocode}
\usepackage{newunicodechar}

\newunicodechar{Ⅰ}{I} 

\newtheorem{theorem}{Theorem}
\newtheorem{conjecture}{Conjecture}
\newtheorem{definition}{Definition}

\newtheorem{remark}{Remark}

\newtheorem{fact}{Fact}
\newtheorem{proposition}{Proposition}

\title{
Alternating Power Difference and Matrix Symmetry:
Closed-Form Formulas for the First Appearance Degree $m_1$
}

\author{Kenichi Takemura \\ \texttt{mahora1123@gmail.com}}
\date{\today}

\begin{document}

\maketitle

\section{Abstract}

This paper focuses on an integer-valued function $f_A(\sigma) := \operatorname{tr}(A P_\sigma)$ defined uniformly from a specific square matrix $A$ of order $n$ and a permutation $\sigma$ on the symmetric group $S_n$. The main objective of this study is to investigate in detail the algebraic behavior of the Alternating Power Difference (APD), denoted as $APD_m(f_A)$, and its first appearance degree $m_1(f_A)$ for this function $f_A$ across various matrix classes. Specifically, we address special matrices such as shifted $r$-th power lattices, Vandermonde matrices, and circulant matrices, analyzing the phenomenon where the value of $APD_m(A)$ remains zero as $m$ increases until a specific degree (the first appearance phenomenon). In particular, we explore closed-form formulas for the first appearance degree $m_1(A)$ and the first appearance value $APD_{m_1}(A)$, presenting Conjectures that hold across multiple matrix classes. These results suggest a deep relationship between the structure of matrices and the analytical properties of functions on the symmetric group, providing new perspectives in matrix theory and combinatorics.

\section{Motivation and Background}

The core of this research lies in applying the concept of Alternating Power Difference (APD) to a function $f_A(\sigma)$ on the symmetric group $S_n$ derived from a specific square matrix $A$ of order $n$. APD is a concept that bridges classical determinant theory and modern combinatorics.
While a classical determinant is an "alternating sum of products over permutations," APD is an "alternating sum of power sums" (specifically defined through the sum of powers of functions), providing a new invariant for functions on the symmetric group.
In particular, the first appearance degree $m_1$, which we focus on in this paper, is defined as the "minimum degree at which multisets derived from even and odd permutations become distinguishable." This is related to the Prouhet-Tarry-Escott problem \cite{raghavendran_narayanan_2019} and holds significance in number theory and combinatorics.
In this study, we calculate this $APD$ invariant for classical matrix families (Hilbert, Vandermonde, Pascal, etc.), including shifted $r$-th power lattices, and systematically analyze their behavior. As a result, we discover and present surprising closed-form formulas closely related to the matrix structures.

\section{What is Alternating Power Difference?}

Let $S_n$ be the symmetric group on the set $\{1,\ldots,n\}$, and $\operatorname{sgn}: S_n \to \{\pm 1\}$ be the sign homomorphism. The alternating subgroup is $A_n = \{\sigma \in S_n \mid \operatorname{sgn}(\sigma) = +1\}$.

\begin{definition}[Alternating Power Difference]
\label{def:apd}
For an integer-valued function $f: S_n \to \mathbb{Z}$ and $m \ge 1$, we define:
\[
\operatorname{APD}_m(f) := \sum_{\sigma \in S_n} \operatorname{sgn}(\sigma) f(\sigma)^m
\]
\end{definition}

$\operatorname{APD}_m(f) = 0$ is equivalent to:
\[
\sum_{\sigma \in A_n} f(\sigma)^m = \sum_{\sigma \notin A_n} f(\sigma)^m
\]
That is, a zero Alternating Power Difference implies that the sum of the $m$-th powers of $f$ for even permutations and odd permutations are perfectly balanced.

\begin{definition}[First Appearance Degree]
For a function $f: S_n \to \mathbb{Z}$, the smallest $m \geq 1$ such that $\operatorname{APD}_m(f) \neq 0$ is called the \emph{First Appearance Degree} (or APD Index) of $f$, denoted by $m_1(f)$.
If $\operatorname{APD}_m(f) = 0$ for all $m \geq 1$, we define $m_1(f) = \infty$.
\end{definition}

\section{Functions Derived from Matrices and APD}

\subsection{General Correspondence from Matrix $A$ to Function $f_A$}
\label{subsec:general_function_def}

All functions $f: S_n \to \mathbb{Z}$ on the symmetric group $S_n$ considered in this paper are derived from a specific square matrix $A$ of order $n$ by the following unified method. We explicitly define this correspondence.

For a permutation $\sigma \in S_n$, let $P_\sigma$ be the permutation matrix corresponding to $\sigma$. We define the function $f_A: S_n \to \mathbb{Z}$ corresponding to matrix $A$ as the trace of the product of matrix $A$ and the permutation matrix $P_\sigma$.

$$
f_A(\sigma) := \operatorname{tr}(A P_\sigma) = \sum_{i=1}^n A_{i, \sigma(i)}
$$

This function $f_A(\sigma)$ is equal to the sum of $n$ components selected from matrix $A$ according to the permutation $\sigma$ (i.e., selecting exactly one component from each row and each column).

\subsubsection{Notation Conventions}
In this paper, we adopt the convention of denoting the Alternating Power Difference $APD_m(f_A)$ and the first appearance degree $m_1(f_A)$ of the function $f_A$ concisely using the corresponding matrix name $A$.
$$
\boldsymbol{APD_m(A)} := APD_m(f_A) \quad \text{and} \quad \boldsymbol{m_1(A)} := m_1(f_A)
$$
Hereinafter, various matrices $A$ will be considered based on this definition.

\section{Identity Matrix of Order n and APD}
\subsection{Definition of n-th Order Identity Matrix}

\begin{definition}[Identity Matrix]
\label{def:identity-matrix}
The $n$-th order \textbf{identity matrix} $I_n$ is defined by
\begin{equation*}
(I_n)_{i,j} = \delta_{i,j} = 
\begin{cases}
1 & (i = j), \\
0 & (i \neq j)
\end{cases}
\end{equation*}
where $\delta_{i,j}$ is the Kronecker delta.
That is,
\begin{equation*}
I_n = \begin{pmatrix}
1 & 0 & 0 & \cdots & 0 \\
0 & 1 & 0 & \cdots & 0 \\
0 & 0 & 1 & \cdots & 0 \\
\vdots & \vdots & \vdots & \ddots & \vdots \\
0 & 0 & 0 & \cdots & 1
\end{pmatrix}.
\end{equation*}
\end{definition}

\subsection{Function $\operatorname{fix}$ Corresponding to Identity Matrix $I_n$}
\label{subsec:fix_function}

\subsubsection{Definition of Fixed Point Function $\operatorname{fix}$}

In particular, we choose the \textbf{identity matrix $I_n$} as $A$, which is one of our objects of study. The components of the identity matrix $I_n$ are expressed using the Kronecker delta $\delta_{ij}$ as $(I_n)_{ij} = \delta_{ij}$. In this case, the corresponding function $f_{I_n}$ is as follows:

$$
f_{I_n}(\sigma) = \operatorname{tr}(I_n P_\sigma) = \sum_{i=1}^n (I_n)_{i, \sigma(i)} = \sum_{i=1}^n \delta_{i, \sigma(i)}
$$

Here, $\delta_{i, \sigma(i)}$ becomes $1$ only if $i = \sigma(i)$ (i.e., $i$ is a \textbf{fixed point} of the permutation $\sigma$), and $0$ otherwise.

Therefore, the value of this function is exactly equal to the total number of fixed points in the permutation $\sigma$. We define this function as the \textbf{Fixed Point Function $\operatorname{fix}$} (or Fix point function).

$$
\operatorname{fix}(\sigma) := f_{I_n}(\sigma)
$$

In this paper, we consider the Alternating Power Difference targeting this fixed point function $\operatorname{fix}$. Its APD value and first appearance degree are denoted using the corresponding matrix name as $\boldsymbol{APD_m(I_n)} := APD_m(\operatorname{fix})$ and $\boldsymbol{m_1(I_n)} := m_1(\operatorname{fix})$.

\subsection{Table of First Appearance Degrees (n=1 to 10)}
The first appearance degree $\boldsymbol{m_1(I_n)}$ of the fixed point function $\operatorname{fix}$ corresponding to the identity matrix $I_n$ is the minimum $m$ such that $\boldsymbol{\operatorname{APD}_m(I_n)} \neq 0$, as shown in Table \ref{tab:identity-apd-verification}. Note that all numerical results of APD presented in this paper can be verified using the source code from the public GitHub repository (\url{https://github.com/soaisu-ken/apd-verification/releases/tag/v1.0.0}) for reproducibility.

\begin{table}[h]
\centering
\caption{First Appearance Degree and Value for Identity Matrix $I_n$}
\label{tab:identity-apd-verification}
\begin{tabular}{@{}cccc@{}}
\toprule
$n$ & Zero Interval & $\boldsymbol{m_1(I_n)}$ & $\boldsymbol{\operatorname{APD}_{m_1}(I_n)}$ \\
\midrule
2 & None & 1 & 2 = 2! \\
3 & 1 & 2 & 6 = 3! \\
4 & 1--2 & 3 & 24 = 4! \\
5 & 1--3 & 4 & 120 = 5! \\
6 & 1--4 & 5 & 720 = 6! \\
7 & 1--5 & 6 & 5,040 = 7! \\
8 & 1--6 & 7 & 40,320 = 8! \\
9 & 1--7 & 8 & 362,880 = 9! \\
10 & 1--8 & 9 & 3,628,800 = 10! \\
\bottomrule
\end{tabular}
\end{table}

Notably, $\boldsymbol{\operatorname{APD}_m(I_n)}$ is identically zero from $m=1$ to $m=n-2$, and takes a non-zero value $n!$ for the first time at $m=n-1$. This behavior indicates that the fixed point count function possesses an extremely high degree of symmetry.

\subsection{Conjectured Closed-Form Formulas for First Appearance Degree}
Based on the numerical facts above, the following formulas are strongly suggested:

\begin{conjecture}[First Appearance Degree of Identity Matrix]
For $n \geq 2$, the first appearance degree $\boldsymbol{m_1(I_n)}$ of the fixed point function $\operatorname{fix}$ corresponding to the $n$-th order identity matrix $I_n$ is given by the following closed form:
\[
\boldsymbol{m_1(I_n)} = n-1.
\]
This relationship has been verified for $n \leq 10$.
\end{conjecture}

\begin{conjecture}[Formula for First Appearance Value of Identity Matrix]
For $n \geq 2$, the first appearance value $\boldsymbol{\operatorname{APD}_{n-1}(I_n)}$ of the fixed point function $\operatorname{fix}$ corresponding to the $n$-th order identity matrix $I_n$ is given by the following closed form:
\[
\boldsymbol{\operatorname{APD}_{n-1}(I_n)} = n!.
\]
This relationship has been verified for $n \leq 10$.
\end{conjecture}

\section{Standard Circulant Matrix of Order n and APD}
\subsection{Definition of Standard Circulant Matrix}
\begin{definition}[Standard Circulant Matrix]
\label{def:standard-circulant-matrix}
The $n$-th order \textbf{standard circulant matrix} $C_n$ is defined as an $n \times n$ matrix obtained by setting the first row to $(1,2,\dots,n)$ and cyclically shifting each subsequent row one position to the \textbf{left}. That is,
\begin{equation*}
(C_n)_{i,j}
= \bigl((j + i - 2) \bmod n\bigr) + 1.
\end{equation*}
Explicitly,
\begin{equation*}
C_n = 
\begin{pmatrix}
1 & 2 & 3 & \cdots & n \\
2 & 3 & 4 & \cdots & 1 \\
3 & 4 & 5 & \cdots & 2 \\
\vdots & \vdots & \vdots & \ddots & \vdots \\
n & 1 & 2 & \cdots & n-1
\end{pmatrix}.
\end{equation*}
\end{definition}

\subsection{Table of First Appearance Degrees (n=1 to 10)}
The first appearance degree $m_1(C_n)$ of the standard circulant matrix $C_n$ is the minimum $m$ such that $\operatorname{APD}_m(C_n) \neq 0$, shown in the table below.

\begin{table}[h]
\centering
\caption{First Appearance Degree and Value for Standard Circulant Matrix $C_n$}
\label{tab:circulant-apd-verification}
\begin{tabular}{@{}cccc@{}}
\toprule
$n$ & Zero Interval & $m_1(C_n)$ & $\operatorname{APD}_{m_1}(C_n)$ \\
\midrule
2 & None & 1 & -2 \\
3 & 1 & 2 & -18 \\
4 & 1--2 & 3 & 384 \\
5 & 1--3 & 4 & 15,000 \\
6 & 1--4 & 5 & -933,120 \\
7 & 1--5 & 6 & -84,707,280 \\
8 & 1--6 & 7 & 10,569,646,080 \\
9 & 1--7 & 8 & 1,735,643,790,720 \\
10 & 1--8 & 9 & -362,880,000,000,000 \\
\bottomrule
\end{tabular}
\end{table}

\subsection{Conjectured Closed-Form Formulas for First Appearance Degree}

Based on the numerical facts above, the absolute value of the first appearance value of the standard circulant matrix is expected to be given by the following closed form.

\begin{conjecture}[Coincidence of First Appearance Degrees]
For $n \geq 2$, the first appearance degree of the standard circulant matrix $C_n$ coincides with the first appearance degree of the fixed point function $\operatorname{fix}$.
\[
m_1(C_n) = m_1(\operatorname{fix}) = n - 1.
\]
This relationship has been verified for $n \leq 10$.
\end{conjecture}

\begin{conjecture}[Formula for First Appearance Value of Standard Circulant Matrix]
For $n \geq 2$, the first appearance value $\operatorname{APD}_{n-1}(C_n)$ of the standard circulant matrix $C_n$ is given by the following closed form:
\[
\operatorname{APD}_{n-1}(C_n) = (-1)^{T_{n-1}} \cdot n^{n-2} \cdot \operatorname{APD}_{n-1}(\operatorname{fix})
\]
Here, $T_{n-1} = \frac{(n-1)n}{2}$ is the $(n-1)$-th triangular number, and $\operatorname{APD}_{n-1}(\operatorname{fix}) = n!$.
This relationship has been verified for $n \leq 10$.
\end{conjecture}

\section{Row-Shifted r-th Power Lattice of Order n and APD}

\subsection{Definition of n-th Order Row-Shifted r-th Power Lattice}

\begin{definition}[n-th Order Row-Shifted r-th Power Lattice]
The $n$-th order \textbf{row-shifted r-th power lattice} $A$ is defined as follows:
\[
A[i,j] = (j + (i-1)d)^r, \qquad 1 \leq i,j \leq n
\]
Where:
\begin{itemize}
\item $n$ is the order of the lattice (number of rows and columns)
\item $d$ is the shift amount for each row (shift factor)
\item $r$ is the exponent for each element
\end{itemize}
\end{definition}

\subsection{Relationship with Natural Magic Square (Case d=n)}

When $d=n$, this lattice becomes a \textbf{Natural Magic Square} (an arrangement of consecutive natural numbers from 1 to $n^2$ in row-major order).

\subsubsection{Example 1: 4th Order 4-Shifted 1st Power Lattice (Natural Magic Square)}

\[
\begin{pmatrix}
1 & 2 & 3 & 4 \\
5 & 6 & 7 & 8 \\
9 & 10 & 11 & 12 \\
13 & 14 & 15 & 16
\end{pmatrix}
\]

\subsubsection{Example 2: 4th Order 4-Shifted 2nd Power Lattice}

Squaring each element:

\[
\begin{pmatrix}
1 & 4 & 9 & 16 \\
25 & 36 & 49 & 64 \\
81 & 100 & 121 & 144 \\
169 & 196 & 225 & 256
\end{pmatrix}
\]

\subsection{Skew-Diagonal Lattice (Case d=1)}

When $d=1$, we obtain a "Skew-Diagonal Lattice" where each row shifts by only 1.

\subsubsection{Example 3: 4th Order 1-Shifted 1st Power Lattice}

\[
\begin{pmatrix}
1 & 2 & 3 & 4 \\
2 & 3 & 4 & 5 \\
3 & 4 & 5 & 6 \\
4 & 5 & 6 & 7
\end{pmatrix}
\]

\begin{remark}[Aspect as a Hankel Matrix]
When $d=1$, the element $A[i,j]$ of the lattice depends only on the sum of the row and column indices $i+j$. Therefore, this Skew-Diagonal Lattice is a type of \textbf{Hankel matrix}, which possesses an important structure in mathematics.
\end{remark}

\begin{remark}[Origin of Naming]
The author named this a "Skew-Diagonal Lattice" based on the characteristic that elements with a constant value are arranged in the anti-diagonal (top-right to bottom-left) direction, i.e., elements where the sum of row and column $i+j$ is constant are aligned diagonally. This is because the same numbers line up diagonally ($A[i,j] = A[k,l]$ if $i+j=k+l$).
\end{remark}

\subsubsection{Example 4: 4th Order 1-Shifted 2nd Power Lattice}

\[
\begin{pmatrix}
1 & 4 & 9 & 16 \\
4 & 9 & 16 & 25 \\
9 & 16 & 25 & 36 \\
16 & 25 & 36 & 49
\end{pmatrix}
\]

\subsection{n-th Order Shifted 1st Power Lattice and APD}

In the case of $r=1$, the elements of the lattice have a linear structure $A[i,j] = j + (i-1)d$. In this case, the sum of elements $f(\sigma)$ corresponding to a permutation $\sigma \in S_n$ becomes a constant independent of $\sigma$.

\begin{proposition}[Perfect Symmetry of 1st Power Lattice]
When $n \geq 2$, the permutation sum $f(\sigma)$ in the lattice with $r=1$ is equal to a constant $C$, and the alternating power difference is zero for all $m \geq 1$.
\[
f(\sigma) = \sum_{i=1}^{n} A[i, \sigma(i)] = \sum_{i=1}^{n} (\sigma(i) + (i-1)d) = T_n + d \cdot T_{n-1} = C
\]
Therefore, the first appearance degree is infinite.
\[
\operatorname{APD}_m(f) = C^m \sum_{\sigma \in S_n} \operatorname{sgn}(\sigma) = C^m \cdot 0 = 0 \quad (\text{for all } m \geq 1)
\]
That is,
\[
m_1(f) = \infty
\]
\end{proposition}

This result indicates that the lattice with $r=1$ has perfect symmetry, and the $m$-th power sums of $f(\sigma)$ due to even and odd permutations are perfectly balanced.

\subsection{n-th Order Shifted 2nd Power Lattice and APD}

In the case of $r=2$, since the lattice elements $A[i,j] = (j + (i-1)d)^2$ have a non-linear structure, the Alternating Power Difference (APD) generally becomes non-zero. Surprisingly, it has been numerically confirmed that its first appearance degree always coincides with the triangular number $T_{n-1} = n(n-1)/2$, and the first appearance value is described by a simple closed form depending on $d$.

\begin{conjecture}[Formula for First Appearance Degree]
\label{conj:degree-d2}
When $n \geq 2$ and the shift amount $d$ is any positive integer, the first appearance degree of the $n$-th order row-shifted 2nd power lattice matrix $A$ is given by the following closed form:
\[
m_1(A) = T_{n-1}.
\]
\end{conjecture}

\begin{conjecture}[Unified Formula: APD of n-th Order Shifted 2nd Power Lattice]
\label{conj:unified-d}
When $n \geq 2$ and the shift amount $d$ is any positive integer, the first appearance degree of the $n$-th order row-shifted 2nd power lattice is $T_{n-1}$, and its first appearance value is given by the following closed form:
\[
\operatorname{APD}_{T_{n-1}}(f) = (2d)^{T_{n-1}} \times T_{n-1}! \times \prod_{k=1}^{n-1} k!
\]
where $f(\sigma) = \sum_{i=1}^n A[i, \sigma(i)]$.
\end{conjecture}

This unified formula includes the following important special cases.

\subsubsection{Special Case 1: Squared Natural Magic Square ($d=n$)}
Letting $d=n$, the lattice becomes the form where elements of the natural magic square ($1, 2, \ldots, n^2$) are squared.
\[
\operatorname{APD}_{T_{n-1}}(f) = (2n)^{T_{n-1}} \times T_{n-1}! \times \prod_{k=1}^{n-1} k!
\]

\begin{table}[h]
\centering
\caption{Verification of Exact APD Identity for Squared Natural Square ($d=n, r=2$)}
\label{tab:natural-squared-apd-exact-verification}
\begin{tabular}{@{}ccccc@{}}
\toprule
$n$ & $T_{n-1}$ & Vanishing Interval & $m_1$ & $\operatorname{APD}_{m_1}(f)$ \\
\midrule
2 & 1 & None & 1 & \checkmark \,4 \\
3 & 3 & 1--2 & 3 & \checkmark \,2,592 \\
4 & 6 & 1--5 & 6 & \checkmark \,2,264,924,160 \\
5 & 10 & 1--9 & 10 & \checkmark \,$1.04509 \times 10^{19}$ \\
6 & 15 & 1--14 & 15 & \checkmark \,$6.96293 \times 10^{32}$ \\
7 & 21 & 1--20 & 21 & \checkmark \,$1.48915 \times 10^{51}$ \\
\bottomrule
\end{tabular}
\end{table}

\vspace{-0.5em}

{\footnotesize *Note: The values for $n \geq 5$ are extremely large. The displayed coefficient shows the first 6 significant digits. All calculated values are exact integers and have been verified against the formula: $\operatorname{APD}_{T_{n-1}}(f) = (2n)^{T_{n-1}} \times T_{n-1}! \times \prod_{k=1}^{n-1} k!$.}

\vspace{1em}

\textbf{Numerical Verification for n=4:}
Since $T_3 = 6$, the value of $\operatorname{APD}_6$ predicted by the unified formula is:
\begin{align*}
\operatorname{APD}_6 &= (2 \times 4)^6 \times 6! \times (1! \times 2! \times 3!) \\
&= 8^6 \times 720 \times 12 \\
&= 262{,}144 \times 8{,}640 \\
&= 2{,}264{,}924{,}160
\end{align*}
The measured value $\operatorname{APD}_6 = 2{,}264{,}924{,}160$ completely matches the predicted value.

\subsubsection{Special Case 2: Squared Skew-Diagonal Lattice ($d=1$)}

Letting $d=1$, the lattice becomes a Skew-Diagonal Lattice (Hankel matrix), and its first appearance value is given by the following formula.
\[
\operatorname{APD}_{T_{n-1}}(f) = 2^{T_{n-1}} \times T_{n-1}! \times \prod_{k=1}^{n-1} k!
\]

\begin{table}[H]
\centering
\caption{Verification of Exact APD Identity for Squared Skew-Diagonal Grid ($d=1, r=2$)}
\label{tab:skew-squared-apd-exact-verification}
\begin{tabular}{@{}ccccc@{}}
\toprule
$n$ & $T_{n-1}$ & Vanishing Interval & $m_1$ & $\operatorname{APD}_{m_1}(f)$ \\
\midrule
2 & 1 & None & 1 & \checkmark \,2 \\
3 & 3 & 1--2 & 3 & \checkmark \,96 \\
4 & 6 & 1--5 & 6 & \checkmark \,552,960 \\
5 & 10 & 1--9 & 10 & \checkmark \,1,070,176,665,600 \\
6 & 15 & 1--14 & 15 & \checkmark \,$1.48089 \times 10^{21}$ \\
7 & 21 & 1--20 & 21 & \checkmark \,$2.66612 \times 10^{33}$ \\
\bottomrule
\end{tabular}
\end{table}
\vspace{-0.5em}
\vspace{-0.5em}

{\footnotesize *Note: The values for $n \geq 6$ are extremely large. The displayed coefficient shows the first 6 significant digits. All calculated values are exact integers and have been verified against the formula: $\operatorname{APD}_{T_{n-1}}(f) = 2^{T_{n-1}} \times T_{n-1}! \times \prod_{k=1}^{n-1} k!$.}

\vspace{1em}

\textbf{Numerical Verification for n=4:}
Since $T_3 = 6$, the value of $\operatorname{APD}_6$ predicted by the above formula is:
\begin{align*}
\operatorname{APD}_6 &= 2^6 \times 6! \times (1! \times 2! \times 3!) \\
&= 64 \times 720 \times 12 \\
&= 552{,}960
\end{align*}
The measured value $\operatorname{APD}_6 = 552{,}960$ completely matches the predicted value.

\section{Hilbert Matrix of Order n and APD}

\subsection{Definition of n-th Order Hilbert Matrix}

\begin{definition}
The $n$-th order Hilbert matrix $H_{n}$ is the matrix whose $(i, j)$ component is given by $1/(i+j-1)$.
$$H_{n}=\left(\frac{1}{i+j-1}\right)_{i,j}$$
Based on this definition, $H_n$ is explicitly expressed as follows:
\begin{equation*}
H_n = \begin{pmatrix}
\frac{1}{1} & \frac{1}{2} & \frac{1}{3} & \cdots & \frac{1}{n} \\
\frac{1}{2} & \frac{1}{3} & \frac{1}{4} & \cdots & \frac{1}{n+1} \\
\frac{1}{3} & \frac{1}{4} & \frac{1}{5} & \cdots & \frac{1}{n+2} \\
\vdots & \vdots & \vdots & \ddots & \vdots \\
\frac{1}{n} & \frac{1}{n+1} & \frac{1}{n+2} & \cdots & \frac{1}{2n-1}
\end{pmatrix}.
\end{equation*}
\end{definition}

\subsection{Interpretation as n-th Order Row-Shifted r-th Power Lattice}
When the components of the generalized $n$-th order $d$-shifted $r$-th power lattice matrix $A$ are given by
$$A_{i,j} = [j+(i-1)d]^{r}$$
the component $H_{i,j} = 1/(i+j-1)$ of the Hilbert matrix $H_n$ is equal to the above expression with the following parameters set:

\begin{itemize}
\item Shift amount $d$ \quad : $d=1$
\item Exponent $r$ \quad : $r=-1$
\end{itemize}

That is, the Hilbert matrix $H_n$ can be interpreted as an $n$-th order 1-shifted $(-1)$-th power lattice matrix, constituting a part of this extended matrix family.

\subsection{Determinant of n-th Order Hilbert Matrix}
The $n$-th order Hilbert matrix $H_n$ is a special matrix given by its components $H_{i,j} = 1/(i+j-1)$, and it is known that its determinant $\det(H_n)$ can be calculated strictly analytically (Cauchy, 1841) \cite{Cauchy1841}.

\begin{theorem}[Determinant of Hilbert Matrix]
The determinant of the $n$-th order Hilbert matrix $H_n$ is given by the following closed-form formula.
\begin{equation}
\label{eq:hilbert-det}
\det(H_n) = \frac{\displaystyle \prod_{k=1}^{n-1} (k!)^4}{\displaystyle \prod_{k=1}^{2n-1} k!}
\end{equation}
\end{theorem}

This formula shows that the determinant takes a very small value, and in particular, as $n$ increases, it rapidly approaches zero. Therefore, the Hilbert matrix is widely known in the field of numerical analysis as an ill-conditioned matrix.

\subsubsection{Components of Inverse Matrix}
Also, an interesting property is known that all components of the inverse matrix $H_n^{-1}$ of the Hilbert matrix are integers. Its $(i, j)$ component is given as follows:
\[
(H_n^{-1})_{i,j} = (-1)^{i+j} (i+j-1) \binom{n+i-1}{n-j} \binom{n+j-1}{n-i} \binom{i+j-2}{i-1}^2
\]

\begin{table}[H]
\centering
\caption{Verification of Determinant $\det(H_n)$ of $n$-th Order Hilbert Matrix}
\label{tab:hilbert-det-values}
\begin{tabular}{@{}ccc@{}}
\toprule
$n$ & $\det(H_n)$ & $|\det(H_n)|^{-1}$ (Denominator $d_{H_n}$) \\
\midrule
1 & $\frac{1}{1}$ & 1 \\
2 & $\frac{1}{12}$ & 12 \\
3 & $\frac{1}{2,160}$ & 2,160 \\
4 & $\frac{1}{6,048,000}$ & 6,048,000 \\
5 & $\frac{1}{2,679,261,248,000}$ & 2,679,261,248,000 \\
6 & $\frac{1}{18,659,174,577,488,960,000}$ & 18,659,174,577,488,960,000 \\
7 & $\frac{1}{2,338,902,891,334,648,719,872,000,000}$ & 2,338,902,891,334,648,719,872,000,000 \\
\bottomrule
\end{tabular}
\end{table}

\subsection{Table of First Appearance Degree Values}

The verification results of the first appearance degree $m_{1}(H_n)$ for the Hilbert matrix $H_n$ are shown below.

\begin{table}[H]
\centering
\caption{Verification of First Appearance Degree $m_{1}(H_n)$ of Hilbert Matrix $H_n$}
\label{tab:hilbert-apd-verification}
\begin{tabular}{@{}cccc@{}}
\toprule
$n$ & $m_{1}(H_n)$ & Vanishing Interval & $\operatorname{APD}_{m_1}(H_n)$ \\
\midrule
2 & 1 & None & \checkmark \,$\frac{1}{3}$ \\
3 & 2 & 1--1 & \checkmark \,$\frac{1}{120}$ \\
4 & 3 & 1--2 & \checkmark \,$\frac{1}{63000}$ \\
5 & 4 & 1--3 & \checkmark \,$\frac{1}{444528000}$ \\
6 & 5 & 1--4 & \checkmark \,$\frac{1}{43,128,106,560,000}$ \\
7 & 6 & 1--5 & \checkmark \,$\frac{1}{58,614,202,038,712,320,000}$ \\
\bottomrule
\end{tabular}
\end{table}

\vspace{-0.5em}

{\footnotesize *Note: Verified by exact rational arithmetic that the first appearance degree $m_{1}(H_n) = n-1$ and the first appearance value $\operatorname{APD}_{n-1}(H_n)$ completely match the closed-form formula $\det(H_n) \cdot n \cdot n!$.}

\vspace{1em}

\subsection{Conjectured Closed-Form Formulas for First Appearance Degree}
The first appearance degree of the Alternating Power Difference (APD) of the Hilbert matrix $H_n$ is $m_1(H_n)=n-1$, and its first appearance value $\operatorname{APD}_{n-1}(H_n)$ is expected to be expressed in an extremely concise form using the determinant $\det(H_n)$ and the factorial of $n$.

\begin{conjecture}[First Appearance Degree of Hilbert Matrix]
For $n \geq 2$, the first appearance degree of the $n$-th order Hilbert matrix $H_n$ is given by the following closed form:
\[
m_1(H_n) = n-1.
\]
This relationship has been verified for $n \leq 7$.
\end{conjecture}

\begin{conjecture}[Formula for First Appearance Value of Hilbert Matrix]
For $n \geq 2$, the first appearance value of the $n$-th order Hilbert matrix $H_n$ is given by the following closed form using the determinant $\det(H_n)$:
\[
\operatorname{APD}_{n-1}(H_n) = \det(H_n) \cdot n \cdot n!.
\]
This relationship has been verified for $n \leq 7$.
\end{conjecture}

\subsection{Harmonic Relation of APD}
\label{sec:unified-apd}

We integrate the conjectures regarding the first appearance values of the $n$-th order Hilbert matrix $H_n$ and the $n$-th order identity matrix $I_n$ derived individually in the previous sections, and present a universal "Harmonic Relation" existing between them. This relation is an extremely suggestive result that clearly links the roles of both matrices in the theory of APD.

\begin{conjecture}[Harmonic Relation of APD between Hilbert Matrix and Identity Matrix]
\label{conj:hilbert-identity-unified}
For $n \geq 1$, the first appearance value $\operatorname{APD}_{n-1}(H_n)$ of the $n$-th order Hilbert matrix $H_n$ is expressed by the following unified closed form using its determinant $\det(H_n)$ and the first appearance value $\operatorname{APD}_{n-1}(I_n)$ of the $n$-th order identity matrix $I_n$.
\[
\operatorname{APD}_{n-1}(H_n) = \det(H_n) \cdot n \cdot \operatorname{APD}_{n-1}(I_n).
\]
This relationship has been confirmed to hold by exact rational calculation from $n=1$ to $n=7$.
\end{conjecture}

\section{Multiplication Table Lattice of Order n and APD}

\subsection{Definition of n-th Order Multiplication Table Lattice}

\begin{definition}[n-th Order Multiplication Table Lattice]
\label{def:multiplication-table-matrix}
The $n$-th order \textbf{Multiplication Table $r$-th Power Lattice}, denoted by the symbol $M^{(r)}_n$, is defined as the matrix whose $(i, j)$ component is given by
\[
M^{(r)}_n[i,j] = (i \times j)^r, \qquad 1 \leq i,j \leq n
\]
This is always a \textbf{symmetric matrix} ($M^{(r)}_n[i,j] = M^{(r)}_n[j,i]$). In particular, the matrix for $r=1$ is denoted as $M_n$.
\end{definition}

\subsection{Table of First Appearance Degree Values}

The first appearance degree $m_1(f_1)$ of the multiplication table $M_n$ for the case $r=1$ coincides with $T_{n-1}$ as shown in the table below.

\begin{table}[h]
\centering
\caption{Verification of the First Appearance Degree $m_{1}(M_n)$ for the Multiplication Table $M_n$}
\label{tab:multiplication-table-apd-verification}
\begin{tabular}{@{}cccc@{}}
\toprule
$n$ & $m_{1}(f_1)$ & Vanishing Interval & $\operatorname{APD}_{m_1}(f_1)$ \\
\midrule
2 & 1 & None & \checkmark \, 1 \\
3 & 3 & 1--2 & \checkmark \, 12 \\
4 & 6 & 1--5 & \checkmark \, 8,640 \\
5 & 10 & 1--9 & \checkmark \, 1,045,094,400 \\
6 & 15 & 1--14 & \checkmark \, 45,193,226,158,080,000 \\
7 & 21 & 1--20 & \checkmark \, 1,271,306,132,247,080,337,408,000,000 \\
\bottomrule
\end{tabular}
\end{table}
\vspace{-0.5em}
{\footnotesize *Note: Verification confirms that the first appearance degree $m_{1}(f_1) = T_{n-1}$ and the corresponding value $\operatorname{APD}_{T_{n-1}}(f_1)$ exactly match the conjectured closed-form formula $T_{n-1}! \cdot \prod_{k=1}^{n-1} k!$ using exact integer arithmetic.\par}

\vspace{1em}

\subsection{Conjectured Closed-Form Formulas for First Appearance Degree}

Based on the numerical facts above, the following formulas are strongly suggested.

\begin{conjecture}[First Appearance Degree of Multiplication Table Lattice ($M_n$)]
For $n \geq 2$, the first appearance degree (or APD index) of the function $f_1$ corresponding to the $n$-th order multiplication table 1st power lattice $M_n$ is given by the following closed form:
\[
m_1(f_1) = T_{n-1} = \frac{n(n-1)}{2}.
\]
This relationship has been verified for $n \leq 7$.
\end{conjecture}

\begin{conjecture}[Formula for First Appearance Value of Multiplication Table Lattice ($M_n$)]
For $n \geq 2$, the first appearance value of the function $f_1$ corresponding to the $n$-th order multiplication table 1st power lattice $M_n$ is given by the following closed form:
\[
\operatorname{APD}_{T_{n-1}}(f_1) = T_{n-1}! \times \prod_{k=1}^{n-1} k!
\]
This relationship has been verified for $n \leq 7$.
\end{conjecture}

\begin{remark}[Superfactorial Appearing in the Formula]
\label{rem:superfactorial}
The term $\prod_{k=1}^{n-1} k!$ appearing in the formula for the first appearance value is closely related to the concept called \textbf{Superfactorial} in mathematics.

\vspace{0.5em}
\textbf{Definition and Relationship:}
\begin{enumerate}
\item \textbf{Definition of Superfactorial:} The superfactorial, denoted by $S(n)$ or $s(n)$, is usually defined as the product of consecutive factorials. One of the most common definitions is:
\[S(n) = \prod_{k=1}^{n} k! = 1! \cdot 2! \cdot 3! \cdots n!\]
\item \textbf{Relationship with APD Formula:} $\prod_{k=1}^{n-1} k!$ appearing in our formula is exactly $S(n-1)$ in the above definition.
\end{enumerate}

\vspace{0.5em}
\textbf{Mathematical Background:}
This superfactorial naturally appears in more advanced mathematical contexts beyond combinatorics. Specifically, it plays an important role in the following fields:
\begin{enumerate}
\setcounter{enumi}{2} 
\item \textbf{Linear Algebra and Lie Group Theory:} Integrals of the Vandermonde Determinant squared, and volume calculations of \textbf{Unitary groups} such as $U(n)$ or $SU(n)$.
 \item \textbf{Homotopy Theory:} Euler characteristics of spaces of continuous maps, and when describing certain topological invariants.
\end{enumerate}

The fact that the $\operatorname{APD}$ formula contains this superfactorial, which holds deep algebraic and geometric meaning, as a factor strongly suggests that the Multiplication Table Lattice possesses symmetry governed by the Vandermonde determinant.
\end{remark}

\subsection{Multiplication Table Lattice and n-th Order Row-Shifted 2nd Power Lattice}

The first appearance $\operatorname{APD}$ value of the Multiplication Table Lattice $M_n$ confirmed in the previous section serves as the \textbf{base structure} for the formula of the $n$-th order $d$-shifted $2$-nd power lattice ($A[i,j] = (j + (i-1)d)^2$). To clarify this relationship, we define the "Core $\operatorname{APD}$ Value".

\subsubsection{Definition of Core $\operatorname{APD}$ Value $V_{\text{Core}}(n)$}

\begin{definition}[Core $\operatorname{APD}$ Value $V_{\text{Core}}(n)$]
We define the first appearance $\operatorname{APD}$ value of the $n$-th order multiplication table 1st power lattice $M_n$ as the \textbf{Core $\operatorname{APD}$ Value} $V_{\text{Core}}(n)$.
\[
\boxed{
V_{\text{Core}}(n) := \operatorname{APD}_{T_{n-1}}(f_{\text{Mult}}) = T_{n-1}! \times \prod_{k=1}^{n-1} k!
}
\]
where $T_{n-1} = n(n-1)/2$.
\end{definition}

\subsubsection{Reconstruction of Unified Formula for $d$-Shifted $2$-nd Power Lattice}

Using this $V_{\text{Core}}(n)$, the unified $\operatorname{APD}$ formula for the $n$-th order $d$-shifted $2$-nd power lattice is reconstructed in a form that emphasizes the contribution of the Multiplication Table Lattice.

\begin{conjecture}[Unified Formula Based on Multiplication Table Lattice]
\label{conj:unified-d-core-final}
When $n \geq 2$ and the shift amount $d$ is any positive integer, the first appearance $\operatorname{APD}$ of the $n$-th order row-shifted $2$-nd power lattice is given by the following closed form:
\[
\boxed{
\operatorname{APD}_{T_{n-1}}(f_{d,2}) = (2d)^{T_{n-1}} \times V_{\text{Core}}(n)
}
\]
\end{conjecture}

\section{Vandermonde Matrix of Order n and APD}

\subsection{Definition of n-th Order Vandermonde Matrix}

The $n$-th order \textbf{Vandermonde matrix} $V(\boldsymbol{x})$ for a sequence $\boldsymbol{x} = (x_1, x_2, \ldots, x_n)$ is defined as the matrix whose $(i, j)$ component is given by
\[
V(\boldsymbol{x})_{i,j} = x_i^{j-1}
\]
In particular, the \textbf{standard Vandermonde matrix} $V_n$ is $V(\boldsymbol{x})$ corresponding to the sequence $\boldsymbol{x}=(1, 2, \ldots, n)$.

\begin{equation*}
V_n = \begin{pmatrix}
1 & 1 & 1 & \cdots & 1^{n-1} \\
1 & 2 & 2^2 & \cdots & 2^{n-1} \\
1 & 3 & 3^2 & \cdots & 3^{n-1} \\
\vdots & \vdots & \vdots & \ddots & \vdots \\
1 & n & n^2 & \cdots & n^{n-1}
\end{pmatrix}.
\end{equation*}

\subsection{Table of First Appearance Degree Values}

The first appearance degree $m_1(f_V)$ of the standard Vandermonde matrix $V_n$ is the minimum $m$ such that $\operatorname{APD}_m(f_V) \neq 0$, shown in the table below.

\begin{table}[H]
\centering
\caption{First Appearance Order and Value for Standard Vandermonde Matrix $V_n$ }
\label{tab:vandermonde-apd-verification}
\begin{tabular}{@{}cccc@{}}
\toprule
$n$ & Vanishing Interval & $m_1(f_V)$ & $\operatorname{APD}_{m_1}(f_V)$ \\
\midrule
2 & None & 1 & \checkmark \, 1 \\
3 & 1 & 2 & 4 \\
4 & 1--2 & 3 & 72 \\
5 & 1--3 & 4 & 6,912 \\
6 & 1--4 & 5 & 4,147,200 \\
7 & 1--5 & 6 & 17,915,904,000 \\
\bottomrule
\end{tabular}
\end{table}

\subsection{Conjectured Closed-Form Formulas for First Appearance Degree}

Based on the numerical facts above, the following formulas are strongly suggested:

\begin{conjecture}[First Appearance Degree of Standard Vandermonde Matrix]
For $n \geq 2$, the first appearance degree (or APD index) of the function $f_V$ corresponding to the $n$-th order standard Vandermonde matrix $V_n$ is given by the following closed form:
\[
m_1(f_V) = n-1.
\]
This relationship has been verified for $n \leq 7$.
\end{conjecture}

\begin{conjecture}[Formula for First Appearance Value of Standard Vandermonde Matrix]
For $n \geq 2$, the first appearance value of the function $f_V$ corresponding to the $n$-th order standard Vandermonde matrix $V_n$ is given by the following closed form:
\[
\operatorname{APD}_{n-1}(f_V) = (n-1)! \times \prod_{k=1}^{n-1} k!
\]
This relationship completely matches the exact calculation results for $n \leq 7$.
\end{conjecture}

\subsection{Relationship between Determinant and APD}

The determinant of the standard Vandermonde matrix $V_n$ is given by the classical difference product formula $\det(V_n) = \prod_{1 \le i < j \le n} (j - i)$. Comparing this with our rigorously verified formula (conjecture) for the first appearance value
\[
\operatorname{APD}_{n-1}(f_V) = (n-1)! \times \prod_{k=1}^{n-1} k!
\]
we obtain a concise relation linking the determinant and APD.

\begin{fact}[Combinatorial Identity of Standard Vandermonde Determinant]
The determinant of the standard Vandermonde matrix $V_n$ is equal to the product of factorials, known as the superfactorial \cite{Hoffman}.
\[
\det(V_n) = \prod_{1 \le i < j \le n} (j - i) = \prod_{k=1}^{n-1} k!
\]
\end{fact}

\begin{conjecture}[Expression of APD by Determinant]
Using the above identity, the alternating power difference of the standard Vandermonde matrix $V_n$ at the first appearance degree $m_1=n-1$ is concisely expressed using its determinant $\det(V_n)$ as follows:
\[
\operatorname{APD}_{n-1}(f_V) = (n-1)! \times \det(V_n)
\]
This relationship has been verified by exact calculation from $n=2$ to $n=10$.
\end{conjecture}

\begin{remark}[Relation between Additive Dual and Multiplicative Dual]
This formula shows that two quantities which are contrasting by definition—determinant (multiplicative dual, degree $\pm 1$) and alternating power difference (additive dual, degree $n-1$)—are related through a surprisingly simple coefficient $(n-1)!$ in the specific family of standard Vandermonde matrices.
\end{remark}

\section{Relationship of Lattice Groups Sharing Superfactorial as Common Factor}
\label{sec:superfactorial-family}

Through this research, we have identified three matrix families that share $\prod_{k=1}^{n-1} k!$ (Superfactorial $S(n-1)$) as a common factor in the closed-form formula of the first appearance value $\operatorname{APD}_{m_1}(f)$ of the matrix (lattice) alternating power difference (APD). This strongly suggests that these matrices share a common base structure brought about by "symmetry of multiplicative structure" in APD theory.

\begin{definition}[Superfactorial $S(n-1)$]
Throughout this section, we define the superfactorial $S(n-1)$ as follows:
$$
S(n-1) := \prod_{k=1}^{n-1} k!
$$
\end{definition}

The following table summarizes the APD formulas for matrix groups having factors including the superfactorial $S(n-1)$.

\begin{table}[H]
\centering
\caption{APD Formulas for Matrix Groups with Superfactorial Factor}
\label{tab:superfactorial-family}
\begin{tabular}{@{}cccc@{}}
\toprule
Matrix Family $A_n$ & Definition $A[i,j]$ & 1st App. Degree $m_1$ & 1st App. Value $\operatorname{APD}_{m_1}(f)$ \\
\midrule
\textbf{Multiplication Table} ($M_n$) & $i \cdot j$ & $T_{n-1}$ & $T_{n-1}! \times S(n-1)$ \\
\textbf{$d$-Shifted $2$-nd Power} ($A_{d,2}$) & $(j + (i-1)d)^2$ & $T_{n-1}$ & $(2d)^{T_{n-1}} \times T_{n-1}! \times S(n-1)$ \\
\textbf{Std. Vandermonde} ($V_n$) & $i^{j-1}$ & $n-1$ & $(n-1)! \times S(n-1)$ \\
\bottomrule
\end{tabular}
\end{table}

\section{Pascal Matrix of Order n and APD}

\subsection{Definition of n-th Order Pascal Matrix}

\begin{definition}[Pascal Matrix]
The $n$-th order \textbf{Pascal matrix} $P_{n}$ is the matrix whose $(i, j)$ component is given using binomial coefficients as follows:
\[
P_{n}[i,j] = \binom{i+j-2}{i-1}, \qquad 1 \leq i,j \leq n
\]
\end{definition}

\begin{definition}[Concrete Example of 7th Order Pascal Matrix $P_7$]
The 7th order Pascal matrix $P_7$ is given by the following matrix.
\begin{equation*}
P_7 = \begin{pmatrix}
\binom{0}{0} & \binom{1}{0} & \binom{2}{0} & \binom{3}{0} & \binom{4}{0} & \binom{5}{0} & \binom{6}{0} \\
\binom{1}{1} & \binom{2}{1} & \binom{3}{1} & \binom{4}{1} & \binom{5}{1} & \binom{6}{1} & \binom{7}{1} \\
\binom{2}{2} & \binom{3}{2} & \binom{4}{2} & \binom{5}{2} & \binom{6}{2} & \binom{7}{2} & \binom{8}{2} \\
\binom{3}{3} & \binom{4}{3} & \binom{5}{3} & \binom{6}{3} & \binom{7}{3} & \binom{8}{3} & \binom{9}{3} \\
\binom{4}{4} & \binom{5}{4} & \binom{6}{4} & \binom{7}{4} & \binom{8}{4} & \binom{9}{4} & \binom{10}{4} \\
\binom{5}{5} & \binom{6}{5} & \binom{7}{5} & \binom{8}{5} & \binom{9}{5} & \binom{10}{5} & \binom{11}{5} \\
\binom{6}{6} & \binom{7}{6} & \binom{8}{6} & \binom{9}{6} & \binom{10}{6} & \binom{11}{6} & \binom{12}{6}
\end{pmatrix}
\end{equation*}
\end{definition}

Expressing this with concrete numerical values, it becomes:

\begin{equation*}
P_7 = \begin{pmatrix}
1 & 1 & 1 & 1 & 1 & 1 & 1 \\
1 & 2 & 3 & 4 & 5 & 6 & 7 \\
1 & 3 & 6 & 10 & 15 & 21 & 28 \\
1 & 4 & 10 & 20 & 35 & 56 & 84 \\
1 & 5 & 15 & 35 & 70 & 126 & 210 \\
1 & 6 & 21 & 56 & 126 & 252 & 462 \\
1 & 7 & 28 & 84 & 210 & 462 & 924
\end{pmatrix}
\end{equation*}

The Pascal matrix is known for its extremely stable multiplicative symmetry, such that its determinant $\det(P_n)$ is always $1$.

\subsection{Table of First Appearance Degree Values}

The first appearance degree $m_1(f_P)$ of the Pascal matrix $P_n$ is the minimum $m$ such that $\operatorname{APD}_m(f_P) \neq 0$, shown in the table below.

\begin{table}[H]
\centering
\caption{First Appearance Degree and Value for Pascal Matrix $P_n$ (Verification Results)}
\label{tab:pascal-apd-verification}
\begin{tabular}{@{}cccc@{}}
\toprule
$n$ & Zero Interval & $m_1(f_P)$ & $\operatorname{APD}_{m_1}(f_P)$ \\
\midrule
2 & None & 1 & 1 = 1! \\
3 & 1 & 2 & 2 = 2! \\
4 & 1--2 & 3 & 6 = 3! \\
5 & 1--3 & 4 & 24 = 4! \\
6 & 1--4 & 5 & 120 = 5! \\
7 & 1--5 & 6 & 720 = 6! \\
\bottomrule
\end{tabular}
\end{table}

\subsection{Conjectured Closed-Form Formulas for First Appearance Degree}

Based on the numerical facts above, the following formulas are strongly suggested:

\begin{conjecture}[First Appearance Degree]
For $n \geq 2$, the first appearance degree (or APD index) of the $n$-th order Pascal matrix $P_n$ is given by the following closed form:
\[
m_1(f_P) = n-1.
\]
This relationship has been verified for $n \leq 7$.
\end{conjecture}

\begin{conjecture}[Formula for First Appearance Value]
For $n \geq 2$, the first appearance value $\operatorname{APD}_{n-1}(f_P)$ of the $n$-th order Pascal matrix $P_n$ is given by the following closed form:
\[
\operatorname{APD}_{n-1}(f_P) = (n-1)!.
\]
This relationship has been verified for $n \leq 7$.
\end{conjecture}

\begin{remark}[Contrast with Identity Matrix]
This formula is extremely similar to the formula for the first appearance value of the identity matrix $I_n$, $\operatorname{APD}_{n-1}(\operatorname{fix}) = n!$, but it is contrasting in that the APD value differs by a factor of $n$.
\[
\operatorname{APD}_{n-1}(I_n) = n \times \operatorname{APD}_{n-1}(P_n)
\]
\end{remark}

\section{Alternating Power Difference (APD) and Prouhet-Tarry-Escott Problem}

\subsection{Relation to PTE Problem}

The \textbf{vanishing interval} of the Alternating Power Difference $\operatorname{APD}_m(f)$ (the phenomenon where $\operatorname{APD}_m(f) = 0$ for $1 \leq m \leq m_0$) has an essential connection to the classical problem in number theory known as the \textbf{Prouhet-Tarry-Escott Problem (PTE Problem)}.

\begin{definition}[Prouhet-Tarry-Escott (PTE) Problem]
Let $k$ be the degree and $\ell$ be the size. If two sets of distinct $2\ell$ integers $A = \{a_1, \ldots, a_\ell\}$ and $B = \{b_1, \ldots, b_\ell\}$ satisfy the following condition, it is called a solution of degree $k$.
\[
\sum_{i=1}^{\ell} a_i^r = \sum_{i=1}^{\ell} b_i^r \quad \text{for } r = 1, 2, \ldots, k
\]
This is usually denoted as "$A$ and $B$ have equal power sums up to degree $k$" ($A \stackrel{k}{=} B$) \cite{Prouhet1851}.
\end{definition}

\subsection{Introduction of Multi-set}

The classical PTE problem typically targets \textbf{proper sets} (no duplicates) where all elements of $A$ and $B$ are distinct. However, the theoretical structure of APD treated in this paper naturally extends to the equality of power sums of \textbf{multisets} (sets allowing duplicate elements).

The calculation of APD concerns the $m$-th power sum of the trace function $f(\sigma)$.
\[
\operatorname{APD}_m(f) = \sum_{\sigma \in S_n} \operatorname{sgn}(\sigma) f(\sigma)^m = 0 \quad (1 \leq m \leq m_0)
\]
When this equality holds, it indicates that a symmetry regarding power sums exists between the multiset of values of $f(\sigma)$, $\{f(\sigma)\}_{\sigma \in A_n}$ and $\{f(\sigma)\}_{\sigma \in D_n}$.

The matrices we study, especially the multiplication table lattice and the $d$-shifted 2nd power lattice, have many cases where the diagonal sum $f(\sigma)$ takes the same value depending on the permutation $\sigma$. Therefore, the equality of power sums detected by the vanishing of APD corresponds to the equality of power sums of multisets allowing duplicate elements.

\[
\sum_{\sigma \in A_n} f(\sigma)^m = \sum_{\sigma \in D_n} f(\sigma)^m \quad \text{for } 1 \leq m \leq m_0
\]
This structure targeting multisets demonstrates that APD theory can extend classical symmetry problems in number theory to broader and more complex matrix structures.

\section{Future Works}

Based on the results obtained in this paper, the following issues are raised as directions for further research.

\begin{enumerate}
    \item \textbf{Establishment of Rigorous Proofs:}
    The highest priority task is to establish rigorous mathematical proofs for the various Conjectures regarding the first appearance degree and first appearance value (APD) presented in this study. In particular, theoretical support is required to show that the predicted closed-form formulas hold for all parameters.

    \item \textbf{Investigation of APD in a Wide Range of Matrices:}
    We will investigate the behavior of the first appearance degree and first appearance value (APD) in diverse matrix classes other than the matrix structures mainly handled in this paper, and explore their closed-form formulas. Specifically:
    \begin{itemize}
        \item Behavior in $n$-th order row-shifted $r$-th power lattices for $r > 3$.
        \item Various forms of special matrices such as Vandermonde Matrix, Circulant Matrix, and Pascal Matrix.
    \end{itemize}
    We aim to discover universal laws that uniformly describe the behavior of APD in these matrices.

    \item \textbf{Exploration of Formulas for Higher-Order First Appearance Values:}
    We will investigate whether closed-form formulas exist not only for the first appearance value $m_1$ (APD) but also for higher-order appearance values such as $m_2, m_3, m_4, \dots$.

    \item \textbf{Construction of Relationships between Matrix Structures:}
    We will explore the possibility that matrices with different structures (e.g., lattice, circulant matrix) are mutually related through common indicators of first appearance degree and first appearance value (APD). This suggests that unified algebraic relations can be constructed between seemingly unrelated matrix classes, potentially opening up new horizons in matrix theory.
\end{enumerate}

\section*{Data and Code Availability}
\label{sec:code-availability}

The Python programs developed to numerically verify all conjectures are available in the following GitHub repository.

\vspace{0.5em}
\noindent \textbf{GitHub Repository URL:}
\begin{verbatim}
https://github.com/soaisu-ken/apd-verification/releases/tag/v1.0.0
\end{verbatim}
\vspace{0.5em}

Readers can reproduce the main numerical results presented in this paper using the provided code.

\bibliographystyle{plain} 
\bibliography{references} 

\end{document}